\documentclass[titlepage,11pt]{article}
\usepackage{latexsym}
\usepackage{amsfonts}
\usepackage{amsmath}
\newcommand\blackslug{\hbox{\hskip 1pt \vrule width 4pt height 8pt depth 1.5pt
        \hskip 1pt}}
\newcommand\bbox{\hfill \quad \blackslug \bigbreak}

\def\ll{,\ldots,}

\title{Near-domination in graphs}
\author{Bruce Reed\\
School of Computer Science, McGill University, Montr\'eal, Canada\\
and the CNRS, France
\\
\\
Alex Scott\thanks{Supported by a
Leverhulme Trust Research Fellowship.}\\
Mathematical Institute, University of Oxford, Oxford OX2 6GG, UK
\\
\\
Paul Seymour\thanks{Supported by ONR grant N00014-14-1-0084 and NSF
grants DMS-1265563 and DMS-1800053.}\\
Princeton University, Princeton, NJ 08544, USA}

\date{May 29 2017; revised \today}

\newtheorem{thm}{}[section]

\newcommand{\Proof}{\noindent{\bf Proof.}\ \ }

\begin{document}
\maketitle
\begin{abstract}
A vertex $u$  of a graph ``$t$-dominates'' a vertex $v$ if there are at most $t$ vertices different from $u,v$ that are adjacent to 
$v$ and not to $u$; and 
a graph is ``$t$-dominating'' if for every pair of distinct vertices, one of them $t$-dominates the other.
Our main result says that if a graph is $t$-dominating, then it is close (in an appropriate sense) to being $0$-dominating.
We also show that an analogous statement for digraphs is false; and
discuss some connections with the Erd\H{o}s-Hajnal conjecture.

\end{abstract}

\section{Introduction}

In this paper, all graphs are finite and have no loops or parallel edges. 
We say a vertex $u$  of a graph {\em $t$-dominates} a vertex $v$ if there are at most $t$ vertices different from $u,v$ 
that are adjacent to $v$ and not to $u$; and a graph is 
{\em $t$-dominating} if for every pair of distinct vertices, one of them $t$-dominates the other.

Graphs that are $0$-dominating are easily understood; they are called ``threshold graphs'' and have several different characterizations,
which we discuss later. But for general fixed $t$, $t$-dominating graphs are not so transparent, and our main result
states that every $t$-dominating graph has bounded ``local difference'' from a $0$-dominating graph. Let us define this.

Let $G, H$ be graphs on the same vertex set. We say that the {\em local difference} between $G,H$ is $d$ if $d$
is the maximum, over all vertices $v$, of 
$$|N_G(v)\setminus N_H(v)| + |N_H(v)\setminus N_G(v)|,$$
where $N_G(v), N_H(v)$ denote the set of neighbours of $v$ in $G,H$ respectively. Thus, if $G,H$ have local difference $d$,
then $H$ can be obtained from $G$ by changing the adjacency of some pairs of vertices, where the changed pairs form a graph
with maximum degree $d$. Local difference is evidently a metric, and could be used, for instance, to describe  
``defective colouring''. A {\em $d$-defective $k$-colouring} of a graph is a partition of its vertex set into $k$ subsets such that
the subgraph induced on each subset has maximum degree at most $d$, and this is the same as saying the graph has local difference
at most $d$ from a $k$-colourable graph.

It is easy to see that if $G$ is $0$-dominating and $G,H$ have local difference at most $d$ then $H$
is $2d$-dominating, and we prove a kind of converse:

\begin{thm}\label{mainthm}
For all integers $t\ge 0$, if $G$ is $t$-dominating then there is a $0$-dominating graph $H$ with the same vertex set 
such that $G,H$ have local difference at most $646t^4$.
\end{thm}

The proof is in three stages: we reduce the problem to ``split graphs'' (graphs with vertex set partitioned into a 
clique and a stable set); we reduce the split graph question to a problem about matrices; and then we solve the matrix problem.
They are carried out in reverse order.

It is natural to ask whether there is an
analogous result for set containment: given a collection $\mathcal F$ of sets such that, 
for all $A,B\in\mathcal F$, either
$A\subset B$ or $B\subset A$, we know immediately that $\mathcal F$ is a chain.  But what if we are only given that
$\min(|A\setminus B|,|B\setminus A|)\le t$?  Is there some $f(t)$ so that we can add/delete at most $f(t)$ 
elements to each
set to obtain a chain?  Equivalently, is there a version of \ref{mainthm} that holds for outneighbourhoods in a
digraph? In section 5, we show that the answer is no.

This research was motivated by the Erd\H{o}s-Hajnal conjecture~\cite{EH}, that for every graph $H$, there exists $\epsilon>0$
such that every graph $G$ not containing $H$ as an induced subgraph has a stable set or clique of cardinality
at least $O(|V(G)|^{\epsilon})$, and we discuss an application of our theorem to that conjecture in the final section.

\section{$0$-domination}
There are several different characterizations of $0$-dominating graphs, as we explain now. A graph $G$ is $0$-dominating if and only if
there are no four distinct vertices $u,u',v,v'$ such that $uu'$ and $vv'$ are edges and the pairs $u,v'$ and $v,u'$ are nonadjacent.
Here $u,v$ might or might not be adjacent, and also $u', v'$ might or might not be adjacent; so this can restated as:
\begin{thm}\label{0dom}
A graph is $0$-dominating if and only if no induced subgraph is isomorphic to $C_4, 2K_2$, or $P_4$.
\end{thm}
($C_4$ denotes the four-vertex cycle graph; $2K_2$ denotes its complement, the graph consisting of two disjoint edges; and $P_4$
denotes the four-vertex path.)
A {\em split graph} is a graph $G$ such that $V(G)$ can be partitioned into a clique and a stable set, and a graph is a split graph
if and only if~\cite{FH} it has no induced subgraph isomorphic to $C_4, 2K_2$ or $C_5$. Thus every $0$-dominating graph is
a split graph. 

A {\em half-graph} is a bipartite graph with bipartition $X,Y$ say, such that $X,Y$ can be ordered as $X=\{x_1\ll x_m\}$
and $Y=\{y_1\ll y_n\}$ with the following property: for all $i,i',j,j'$ with $1\le i\le i'\le m$ and $1\le j\le j'\le n$, if $x_iy_j$
is an edge then $x_{i'}y_{j'}$ is an edge. It is easy and well-known that a bipartite graph is a half-graph if and only if it has
no induced subgraph isomorphic to $2K_2$. Let us say a {\em split half-graph} is a graph obtained from a half-graph
with bipartition $X,Y$ by adding edges to make $X$ a clique.

Let $G$ be $0$-dominating, and let $V(G)$ be the disjoint union of a clique $X$ and a stable set $Y$. Let $B$ be
the bipartite
graph with bipartition $X,Y$ formed by the edges of $G$ between $X$ and $Y$. Then $B$ has no induced subgraph $2K_2$, since
$G$ has no $P_4$, and it follows
that $B$ is a half-graph, and so $G$ is a split half-graph. Since split half-graphs are $0$-dominating, we have shown that:
\begin{thm}\label{0dom1}
$G$ is $0$-dominating if and only if $G$ is a split half-graph.
\end{thm}

Every nonnull split half-graph has either a vertex of degree zero or a vertex adjacent to all other vertices, and deleting this vertex
gives another split half-graph. A graph $G$ is a {\em threshold graph} if it can be built starting from the null graph by
the two operations of adding an isolated vertex and adding a vertex adjacent to all current vertices; and since threshold graphs
are $0$-dominating, we have a third characterization:
\begin{thm}\label{0dom2}
$G$ is $0$-dominating if and only if $G$ is a threshold graph.
\end{thm}

The characterization \ref{0dom1} is the most useful for our purposes.

\section{A matrix problem}

Let us denote the set of all pairs $(i,j)\;(1\le i\le m, 1\le j\le n)$
by $[m]\times [n]$.
A subset $F$ of $[m]\times [n]$ is {\em up-closed} if for all $i,i', j,j'$ with $1\le i\le i'\le m$
and $1\le j\le j'\le n$, if $(i,j)\in F$ then $(i',j')\in F$; and $F$ is {\em down-closed} if $([m]\times [n])\setminus F$ is up-closed.
The union of two up-closed subsets is up-closed, and so
every subset of $[m]\times [n]$ has a unique maximal up-closed subset.

Now let $A=(a_{ij}:(i,j)\in [m]\times [n])$ be a 0/1 matrix. Its {\em support} is the set of all pairs $(i,j)\in [m]\times [n]$ with $a_{ij} = 1$. 
We say that $A$ is {\em monotone} if its support is up-closed.
For $i,i'\in \{1\ll m\}$ we say that $i,i'$ are {\em row-inclusive} in $A$ if one of the sets $\{j\in \{1\ll n\}:a_{ij}\ne \emptyset\}$,
$\{j\in \{1\ll n\}:a_{i'j}\ne \emptyset\}$ is a subset of the other (and we define {\em column-inclusive}
similarly).
We say that the matrix $A$ is {\em inclusive} if all $i,i'\in \{1\ll m\}$ are row-inclusive in $A$ (and therefore all $j,j'\in \{1\ll n\}$
are column-inclusive in $A$).
It follows that $A$ is inclusive if and only its rows and columns
can be reordered to make a matrix that is monotone.

If $t\ge 0$ is an integer, we say that $A$ is {\em $t$-restricted} if 
\begin{itemize}
\item for all $i,i'\in \{1\ll m\}$ with $i<i'$, there are at most $t$ values of $j\in \{1\ll n\}$ such that $a_{ij} = 1$ and $a_{i'j} = 0$; and
\item for all $j,j'\in \{1\ll n\}$ with $j<j'$, there are at most $t$ values of $i\in \{1\ll m\}$ such that $a_{ij} = 1$ and $a_{ij'} = 0
$.
\end{itemize}
Thus a matrix $A$ is $0$-restricted if and only if it is monotone, that is, its support is up-closed.
In order to prove \ref{mainthm} we need to prove something similar for $t$-restricted matrices; but before that we handle a special case.

Let $A,B$ be two 0/1 matrices both indexed by $[m]\times [n]$. The {\em local difference} between $A,B$ is the maximum, over all rows and columns,
of the number of terms in that row or column in which $A,B$ differ.
One might hope that for all $t$, every  $t$-restricted 0/1 matrix $A$
has bounded local difference from some monotone matrix; but that is false, even for $t=1$.
For instance,
with $m=1$ and $n$ large, a $1\times n$ matrix $A$ with entries $(1\ll 1,0\ll 0)$ (with $n/2$ ones and $n/2$ zeroes)
is $1$-restricted, and yet has arbitrarily large local difference from
every monotone matrix. This is perhaps unfair in that being 1-restricted has no content for a $1\times n$ matrix,
but we could pad it by adding more rows; say $n/2$ rows of all zeroes, then the given row, and then $n/2$ rows of all ones, making
an $(n+1)\times n$ matrix, which is also a counterexample. One might try assuming in addition that the rows are in increasing order of
row-sum, and the same for columns (which will be the case when we apply these results to our graph problem), but a similar
counterexample still can be made, as follows. Take the $(n+1)\times n$ matrix just described, and change the $n/2$ entries $a_{i,\frac{n}{2}+i}\;(1\le i\le n/2)$ to ones
and the $n/2$ entries $a_{\frac{n}{2}+i+1,i}\;(1\le i\le n/2)$ to zeroes.

Nevertheless, something like this is true; replace ``monotone'' by ``inclusive''. We will prove:
\begin{thm}\label{matrix}
For all integers $t\ge 0$ and every $t$-restricted {\rm 0/1} matrix $A$,
there is an inclusive matrix $B$, such that the local difference between $A,B$
is at most $644t^4$.
\end{thm}
This means that we can change a bounded number of entries in every row and column and then reorder rows and columns to
get a monotone matrix.
But before we prove this, we handle a special case, matrices of bounded ``breadth'', and next we define this.
Let $A=(a_{ij}:(i,j)\in E)$ be a 0/1 matrix, with support $E_1$, and let $E_0= ([m]\times [n])\setminus E_1$. 
Let $X$ be the maximal down-closed subset of $E_0$, and $Y$ the maximal up-closed subset of $E_1$. ($X$ is unique, since the union
of two down-closed sets is also down-closed, and similarly $Y$ is unique.) Thus $X\cap Y=\emptyset$.
A {\em diagonal} for $A$ means a subset $F$ of $[m]\times [n]$ such that for some integer $c$, $F$ is the set of all pairs $(i,j)\in [m]\times [n]$
with $j=i+c$ and $(i,j)\notin X\cup Y$. Since $X$ is down-closed and $Y$ is up-closed, it follows 
that every diagonal $F$ is an interval in the sense that if 
$(i,i+c), (i'',i''+c)\in F$ where $i''>i$ then $(i',i'+c)\in F$ for all $i'$ with $i\le i'\le i''$. We call the maximum cardinality
of all diagonals the {\em breadth} of $A$. We first show:

\begin{thm}\label{matrix1}
For all integers $t,w\ge 0$ with $2w\ge t+1$, and every $t$-restricted {\rm 0/1} matrix $A$ with
breadth at most $w$,
there is an inclusive matrix $B$ such that the local difference between $A,B$ is at most $2(t+w)w^3$.
\end{thm}
\Proof
Let $A$ be a $t$-restricted 0/1 matrix indexed by $[m]\times [n]$, with breadth at most $w$, let its support be $E_1$, 
let $E_0=([m]\times [n])\setminus E_1$, and let $X,Y$ be 
the maximal down-closed
subset of $E_0$ and up-closed subset of $E_1$ respectively. Let $Z=([m]\times [n])\setminus (X\cup Y)$.
We will give rules to change some of the entries $a_{ij}$ for $(i,j)\in Z$. We need to satisfy two conditions:
\begin{itemize}
\item after the changes, all pairs $i,i'$ will be row-inclusive; and
\item for each row or column, at most $(2t+w)w^3$ entries in that row or column will be changed.
\end{itemize}
Our first task is to give the rules, but that needs a number of definitions.
For $1\le j\le n$, the $j$th {\em column} means the set $\{(i,j)\;:1\le i\le m\}$ for some $j\in \{1\ll n\}$,
and a {\em row} is defined similarly. A {\em post} is the nonempty intersection of a column with $Z$, and a {\em beam}
is the nonempty intersection of a row with $Z$. Every member of $Z$ belongs to a unique post and a unique beam,
and how we change the corresponding entry of $A$ depends on the types of this post and beam. (Posts and beams again are intervals,
in the natural sense.) 

Let $P,P'$ be posts, where $P, P'$ are the intersection of the $j$th column with $Z$ and the $j'$th column with $Z$, respectively.
They are {\em parallel} if 
\begin{itemize}
\item for all $i\in \{1\ll m\}$, $(i,j)\in Z$ if and only if $(i,j')\in Z$; and
\item for all $i\in \{1\ll m\}$, if $(i,j)\in Z$ then $a_{ij} = a_{ij'}$.
\end{itemize}
Thus two parallel posts involve the same rows, and have the same entries in those rows.
The {\em multiplicity} of a post $P$ is the number of posts that are parallel to $P$
(counting $P$ itself).
We define parallelness and multiplicity for beams similarly.

Let $(i,j)\in Z$. We associate four integers with $(i,j)$:
\begin{itemize}
\item $p^-(i,j)$, the number of $i'<i$ such that $(i',j)\in Z$;
\item $q^-(i,j)$, the number of $j'<j$ such that $(i,j')\in Z$;
\item $p^+(i,j)$, the number of $i'>i$ such that $(i',j)\in Z$;
\item $q^+(i,j)$, the number of $j'>j$ such that $(i,j')\in Z$.
\end{itemize}
We observe that the post containing $(i,j)$ has cardinality $p^-(i,j)+p^+(i,j)+1$, and a similar statement holds for the beam.
\\
\\
(1) {\em For each $(i,j)\in Z$, $\min(p^-(i,j), q^+(i,j))<w$, and $\min(q^-(i,j), p^+(i, j)) <w$.}
\\
\\
Suppose that $p^-(i,j), q^+(i,j)\ge w$. It follows that $(i-w,j), (i,j+w)\in Z$.
Since $X$ is down-closed and $(i-w,j)\notin X$, it follows that for $0\le h\le w$,
$(i-w+h,j+h)\notin X$; and similarly, since $(i,j+w)\notin Y$ it follows that for $0\le h\le w$,
$(i-w+h,j+h)\notin Y$. Consequently the set $\{(i-w+h,j+h)\::(0\le h\le w)\}$ is a subset of a diagonal. 
But it has cardinality $w+1$,
a contradiction. Consequently $\min(p^-(i,j), q^+(i,j))<w$, and similarly $\min(p^+(i,j), q^-(i, j)) <w$.
This proves (1).
\\
\\
(2) {\em For every post $P$, if $(i,j)\in P$ then the multiplicity of $P$ is at most $q^-(i,j)+q^+(i,j)+1$.}
\\
\\
The beam containing $(i,j)$ has cardinality $q^-(i,j)+q^+(i,j)+1$; but it intersects all posts parallel to $P$, and 
so its cardinality is at least the multiplicity of $P$.
This proves (2).

\bigskip
Now we can give the rules. For each $(i,j)\in [m]\times [n]$ define $b_{ij}$ as follows.
If $(i,j)\notin Z$ then $b_{ij} = a_{ij}$, so we may assume that $(i,j)\in Z$; let $(i,j)$ belong to a post $P$ and a beam $Q$.
\begin{itemize}
\item If both $p^-(i,j), q^-(i,j)\ge w$ then $b_{ij} = 1$;
\item if both $p^+(i,j), q^+(i,j)\ge w$ then $b_{ij} = 0$;
\item if both $p^-(i,j), p^+(i,j)<w$ and $P$ has multiplicity at least $2w$ then $b_{ij} = a_{ij}$;
\item if both $q^-(i,j), q^+(i,j)<w$ and $Q$ has multiplicity at least $2w$ then $b_{ij} = a_{ij}$;
\item if both $p^-(i,j), p^+(i,j)<w$ and $P$ has multiplicity less than $2w$ then $b_{ij} = 0$;
\item if both $q^-(i,j), q^+(i,j)<w$ and $Q$ has multiplicity less than $2w$ then $b_{ij} = 0$.
\end{itemize}
We claim the rules are consistent; for let $(i,j)\in Z$. By (1), only one of the first two rules applies to $(i,j)$,
and if one of the first two applies then none of the other rules apply. If the third rule applies to $(i,j)$, then the fifth does not;
and since $P$ has multiplicity 
at least $2w$, it follows by (2) that $q^-(i,j)+q^+(i,j)+1\ge 2w$, so at least one of $q^-(i,j),q^+(i,j)\ge w$ and 
consequently the fourth and sixth rules
do not apply. Similarly if the fourth rule applies then the fifth and sixth do not. Finally, both the fifth and sixth may
apply simultaneously, but they assign the same value to $b_{ij}$. Thus the rules are consistent.
Furthermore, we observe that every $(i,j)\in Z$ falls under one of the rules, by (1), and so the matrix $B=(b_{ij})$ is well-defined.

Now we must check the two bullets given at the start of this proof.
We say $B$ is {\em increasing in row $i$} if $b_{ij}\le b_{ij'}$ for all $j,j'$ with $1\le j<j'\le n$,
and we define {\em increasing in column $j$} similarly.
\\
\\
(3) {\em Let $P$ be a post, a subset of the $j$th column. If $|P|<2w$ then $B$ is increasing in column~$j$.}
\\
\\
Suppose that there exist $i,i'$ with $1\le i<i'\le m$, such that $b_{ij} = 1$ and $b_{i'j} = 0$. It follows that
$(i,j)\notin X$, and so $(i',j)\notin X$, since $X$ is down-closed; and similarly $(i,j), (i',j)\notin Y$, and consequently
$(i,j), (i',j)\in Z$. Let $Q, Q'$ be the beams containing $(i,j), (i',j)$ respectively. Since $|P|<2w$, both
$Q,Q'$ have multiplicity less than $2w$. Since $b_{ij} = 1$, not both $q^-(i,j), q^+(i,j)<w$ from the sixth rule. 
If $p^+(i,j)\ge w$,
then by (1) $q^-(i,j)<w$, so $q^+(i,j)\ge w$ and $b_{ij}= 0$, a contradiction. Thus $p^+(i,j)<w$. 

Suppose that also $p^-(i,j)<w$.
Since $b_{ij} = 1$ it follows that $P$ has multiplicity at least $2w$, and $a_{ij} = b_{ij} = 1$. If $a_{i'j}=0$, 
then since $P$ has multiplicity at least $2w\ge t+1$,
there are $t+1$ values of $j'\in \{1\ll n\}$ such that
$a_{ij'}=1$ and $a_{i'j'}=0$, contradicting that $A$ is $t$-restricted. Thus $a_{i'j} = 1$, and in particular
$a_{i'j}\ne b_{i'j}$. Since $P$ has multiplicity at least $2w$ it follows that one of $p^-(i',j), p^+(i',j)\ge w$, and hence
$p^-(i',j)\ge w$, since $p^+(i',j)<p^+(i,j)<w$. By (1), $q^+(i',j)<w$. Since $P$ has multiplicity at least $2w$,
(2) implies that $q^-(i',j)+q^+(i',j)+1\ge 2w$, and so $q^-(i',j) \ge w$. This contradicts the first rule,
since $b_{i'j} = 0$. This proves that $p^-(i,j)\ge w$. 

We have seen that 
not both $q^-(i,j), q^+(i,j)<w$, and so from (1), $q^-(i,j)\ge w$. But $q^-(i',j)\ge q^-(i,j)$ since
$X$ is down-closed, and so $q^-(i',j)\ge w$; and also $p^-(i',j)>p^-(i,j)\ge w$, and yet $b_{i'j}=0$,
contrary to the rules. This proves (3).
\\
\\
(4) {\em Let $1\le i<i'\le m$; then $i,i'$ are row-inclusive in $B$.}
\\
\\
Suppose not, then in particular, there exists $j\in \{1\ll n\}$ such that $b_{ij} = 1$ and $b_{i'j}=0$. 
Now $(i,j)\notin X$ since $b_{ij} = 1$, and since $X$ is down-closed it follows that $(i',j)\notin X$. Also
$(i',j)\notin Y$ since $b_{i'j} = 0$, and so $(i,j)\notin Y$ since $Y$ is up-closed. 
Consequently $(i,j), (i',j)\in Z$. Let $P$
be the post containing them both, and let $Q, Q'$ be the beams containing $(i,j), (i',j)$ respectively. 
Since $B$ is not increasing in the $j$th column, $|P|\ge 2w$ by (3). If $B$ is increasing in the $i$th row, and $B$
is increasing in the $i'$th row, then (since $B$ is $0/1$-valued) $i,i'$ are row-inclusive in $B$, a contradiction. So 
$B$ is not increasing in one of the $i$th row or the $i'$th row, and so one of $|Q|,|Q'|\ge 2w$ by (3).

Suppose that $|Q|\ge 2w$. Then one of $q^-(i,j), q^+(i,j)\ge w$, and one of $p^-(i,j), p^+(i,j)\ge w$,
and since $b_{ij}=1$, it follows from (1) and the rules that $p^-(i,j), q^-(i,j)\ge w$.
But $q^-(i',j)\ge q^-(i,j)$ since
$X$ is down-closed, and so $q^-(i',j)\ge w$; and also $p^-(i',j)>p^-(i,j)\ge w$, and yet $b_{i'j}=0$,
contrary to the first rule. Thus $|Q|<2w$, and so $|Q'|\ge 2w$. 

Hence one of $q^-(i',j), q^+(i',j)\ge w$, and one of $p^-(i',j), p^+(i',j)\ge w$,
and since $b_{i'j}=0$, it follows from (1) and the rules that $p^+(i',j), q^+(i',j)\ge w$.
But $q^+(i,j)\ge q^+(i',j)$,
and so $q^+(i,j)\ge w$; and also $p^+(i,j)>p^+(i',j)\ge w$, and yet $b_{ij}=1$,
contrary to the rules. This proves (4).
\\
\\
(5) {\em The local difference between $A,B$ is at most $2(t+w)w^3$.}
\\
\\
Let $1\le j\le n$; by the symmetry between axes, it suffices to show that there are at most $(2t+w)w^3$ values 
of $i$ such that $a_{ij}\ne b_{ij}$. For every such value of $i$, it follows that $(i,j)\in Z$; so we may
assume that there is a post $P$ included in the $j$th column with $|P|> (2t+w)w^3$, 
and choose $i_1,i_2$ with $1\le i_1\le i_2\le m$ such that 
\begin{itemize}
\item for $1\le i<i_1$, $(i,j)\in X$;
\item for $i_1\le i\le i_2$, $(i,j)\in P$; and
\item for $i_2<i\le m$, $(i,j)\in Y$.
\end{itemize}
Let $i_1\le i\le i_2$, and suppose that $a_{ij} = 0$ and $b_{ij} = 1$. From the rules, it follows that $p^-(i,j),q^-(i,j)\ge w$;
and so from (1), $p^+(i,j)<w$. Consequently $i_2-w+1\le i\le i_2$, and so there are at most $w$ such values of $i$.
Now let $i_1\le i\le i_2$, and suppose that $a_{ij} = 1$ and $b_{ij} = 0$. From the rules, it follows that either 
\begin{itemize}
\item $p^+(i,j),q^+(i,j)\ge w$, or 
\item both $p^-(i,j), p^+(i,j)<w$ and $P$ has multiplicity less than $2w$, or 
\item both $q^-(i,j), q^+(i,j)<w$ and $Q_i$ has multiplicity
less than $2w$, where $Q_i$ is the beam containing $(i,j)$. 
\end{itemize}
At most $w$ values of $i$ satisfy the first bullet, as before. If $i$ satisfies the second bullet then $|P|\le 2w-1$, a contradiction
since $|P|> (2t+w)w^3$;
so we may assume that $i$ satisfies the third bullet. There are therefore only $w^2$ possibilities for the pair
$q^-(i,j), q^+(i,j)$. Let $h,k\ge 0$ with $h,k<w$, and let $I(h,k)$ be the set of all $i$ that satisfy
$i_1\le i\le i_2$, and $a_{ij}=1$, and $b_{ij}=0$, and $q^-(i,j)=h$ and $q^+(i,j) = k$, and $Q_i$ has multiplicity less than $2w$.
We need to bound $|I(h,k)|$.
We can evidently get an exponential bound, since there are only $2^{2w+1}$ possibilities for the entries of $Q_i$, and 
only $2w-1$ values of $i\in I(h,k)$ in which all the entries are the same; but we can do better.

If $i\in I(h,k)$, $A$ may or may not be increasing in the $i$th row. Let $I_1$ be the set of $i\in I(h,k)$ such that 
$A$ is increasing in the $i$th row, and $I_2=I(h,k)\setminus I_1$. For each $i\in I_1$, there are only $h+1$ possibilities for
the entries of $Q_i$, and at most $2w-1$ distinct $i\in I_1$ with the same set of entries, and so
$|I_1|\le (h+1)(2w-1)$.
If $i\in I_2$ then there exists $j'$ with 
$j-h\le j'\le j+k-1$ such that $a_{ij'} = 1$ and $a_{i, j'+1} = 0$. 
There are only $h+k$ choices for $j'$, and for each $j'$ there are at most $t$
values of $i\in I_2$ such that $a_{ij'} = 1$ and $a_{i, j'+1} = 0$, since $A$ is $t$-restricted; and so $|I_2|\le t(h+k)$. 
Consequently $|I(h,k)|\le (h+1)(2w-1)+t(h+k)$. Since $h,k< w$, it follows that $|I(h,k)|\le w(2w-1)+2(w-1)t$.

In summary, at most $w$ values of $i$ satisfy the first bullet above; none satisfy the second bullet; 
and (since there are
at most $w^2$ choices for $h,k$), at most
$(w(2w-1)+2(w-1)t)w^2$ satisfy the third. Consequently there are at most $(w(2w-1)+2(w-1)t)w^2 + w$ values of $i$
such that $a_{ij}=1$ and $b_{ij}=0$. As we saw, there are at most $w$ values of $i$ such that $a_{ij} = 0$ and $b_{ij}=1$; so there 
are at most $(w(2w-1)+2(w-1)t)w^2 + 2w$ values of $i$ such that $a_{ij}\ne b_{ij}$. This proves (5).

\bigskip

From (4) and (5), this proves \ref{matrix1}.~\bbox

\bigskip

Now we eliminate the hypothesis about bounded breadth, by means of the following.
\begin{thm}\label{matrix2}
For all integers $t\ge 0$ and every $t$-restricted {\rm 0/1} matrix $A$,
there is a $t$-restricted {\rm 0/1} matrix $B$ with breadth at most $4t$, such that the local difference between $A,B$
is at most $4t$.
\end{thm}
\Proof
Let $A$ be a 0/1 matrix indexed by $[m]\times [n]$, let its support be $E_1$,
let $E_0=([m]\times [n])\setminus E_1$, and let $X,Y$ be
the maximal down-closed
subset of $E_0$ and up-closed subset of $E_1$ respectively. Let $Z=([m]\times [n])\setminus (X\cup Y)$.
Let $P_1$ be the set of all $(i,j)\in [m]\times [n]$ such that                    
there exist  at least $2t$ values of $i'< i$ with $(i',j)\in E_1$, and let 
$Q_1$ be the set of all $(i,j)$ such that there are at least $2t$ values of $j'< j$ with $(i,j')\in E_1$.
Let $P_0$ be the set of all $(i,j)$ 
such that there exist  at least $2t$ values of $i'> i$ with $(i',j)\in E_0$, and let $Q_0$ be 
the set of all $(i,j)$ such that there are at least $2t$ values of $j'> j$ with $(i,j')\in E_0$.
\\
\\
(1) {\em There do not exist $(i_1,j_1)\in P_1$ and $(i_0,j_0)\in Q_0$ such that $i_0\ge i_1$ and $j_0\ge j_1$. Also,
there do not exist $(i_1,j_1)\in Q_1$ and $(i_0,j_0)\in P_0$ such that $i_0\ge i_1$ and $j_0\ge j_1$.}
\\
\\
Suppose that such $(i_1,j_1),(i_0,j_0)$ exist. From the symmetry between axes, we may assume that 
$(i_1,j_1)\in P_1$, and hence $(i_0,j_0)\in Q_0$. Consequently $(i_0,j_1)\in P_1\cap Q_0$.
Choose $I_1\subseteq \{1\ll i_0-1\}$ with cardinality $2t$ 
such that $(i,j_1)\in E_1$ for each $i\in I_1$, and choose $J_0\subseteq \{j_1+1\ll n\}$
with cardinality $2t$ such that $(i_0,j)\in E_0$ for each $j\in J_0$.
From the symmetry between zeroes and ones, we may assume without loss of generality that $(i_0,j_1)\in E_0$.
If $E_1$ contains at least half of the pairs $(i,j)$ with $i\in I_1$ and $j\in J_0$, then there exists
$i\in I_1$ such that $(i,j)\in E_1$ for at least $t$ values of $j\in J_0$, and hence for at least
$t+1$ values of $j\in J_0\cup \{j_0\}$; and since 
$(i,j_0)\in E_0$ for each such $j$, this contradicts that $A$ is $t$-restricted. On the other hand, if 
$E_0$ contains more than half of the pairs $(i,j)$ with $i\in I_1$ and $j\in J_0$, then there exists $j\in J_0$
such that $(i,j)\in E_0$ for at least $t+1$ values of $i\in I_1$;  and since
$(i_1,j)\in E_1$ for each such $i$, this also contradicts that $A$ is $t$-restricted. This proves~(1).
\\
\\
(2) {\em There do not exist $(i_1,j_1)\in P_1\cup Q_1$ and $(i_0,j_0)\in P_0\cup Q_0$ such that $i_0\ge i_1+4t$ and $j_0\ge j_1+4t$.}
\\
\\
Suppose that such $(i_1,j_1),(i_0,j_0)$ exist. From the symmetry between axes, we may assume that 
$(i_1,j_1)\in P_1$. By (1) it follows that $(i_0,j_0)\notin Q_0$, and so $(i_0,j_0)\in P_0$.
Choose $I_1\subseteq \{1\ll i_1-1\}$ with cardinality $2t$
such that $(i,j_1)\in E_1$ for each $i\in I_1$, and choose $I_0\subseteq \{i_0+1\ll m\}$
with cardinality $2t$ such that $(i,j_0)\in E_0$ for each $i\in I_0$. If $E_0$ contains more than half of the pairs $(i,j)$
with $i\in I_1$ and $j_1<j< j_0$, then there exists $j$ with $j_1< j< j_0$ such that $(i,j)\in E_0$ for at least $t+1$
values of $i\in I_1$, contradicting that $A$ is $t$-restricted. If $E_1$ contains at least half
of the pairs $(i,j)$ with $i\in I_1$ and $j_1<j<j_0$, then there exists $i\in I_1$ such that $(i,j)\in E_1$
for at least half of the values of $j$ with $j_1< j< j_0$, and since $(i,j_1)\in E_1$, it follows that $(i,j_0)\in Q_1$, contrary to (1)
since $(i_0,j_0)\in P_0$. This proves (2).

\bigskip

Construct a matrix $B = (b_{ij})$ as follows. Let $(i,j)\in [m]\times [n]$. 
\begin{itemize}
\item If $(i,j)\in Z$ and there is no $(i_0,j_0)\in P_0\cup Q_0$ such that $i_0\ge i$ and $j_0\ge j$ then $b_{ij} = 1$;
\item if $(i,j)\in Z$ and there exists $(i_0,j_0)\in P_0\cup Q_0$ such that $i_0\ge i$ and $j_0\ge j$, 
and there is no $(i_1,j_1)\in P_1\cup Q_1$ such that 
$i_1\le i$ and $j_1\le j$, then $b_{ij} = 0$;
\item if neither of these applies then $b_{ij} = a_{ij}$.
\end{itemize}
We need to show that $B$ is $t$-restricted, $B$ has breadth at most $4t$, and the local difference between $A,B$ is at most $4t$. 
Let the support of $B$ be $F_1$,
let $F_0=([m]\times [n])\setminus F_1$, and let $X',Y'$ be
the maximal down-closed
subset of $F_0$ and up-closed subset of $F_1$ respectively.
\\
\\
(3) {\em For all $(i,j)\in [m]\times[n]$, if there is no $(i_0,j_0)\in P_0\cup Q_0$ such that $i_0\ge i$ and $j_0\ge j$ 
then $(i,j)\in Y'$. Also, if there is no $(i_1,j_1)\in P_1\cup Q_1$ such that
$i_1\le i$ and $j_1\le j$, then $(i,j)\in X'$.}
\\
\\
Suppose there is no $(i_0,j_0)\in P_0\cup Q_0$ such that $i_0\ge i$ and 
$j_0\ge j$. Let $F$ be the set of all $(i',j')\in [m]\times[n]$ such that $i\le i'$ and $j\le j'$; then for every such
pair, $(i',j')\in Z\cup Y$, and there is no $(i_0,j_0)\in P_0\cup Q_0$ such that $i_0\ge i'$ and 
$j_1\ge j'$, and so $(i',j')\in F_1$. Hence $F$ is an up-closed subset of $F_1$, 
and since $Y'$ is the unique maximal  up-closed
subset of $F_1$, it follows that $F\subseteq Y'$, and in particular $(i,j)\in Y'$ as required. The second statement follows
similarly. This proves (3).
\\
\\
(4) {\em $B$ is $t$-restricted.}
\\
\\
Suppose that there exist $i,i'$ with $1\le i<i'\le m$ such that $b_{ij} = 1$ and $b_{i'j} = 0$ for $t+1$ values of $j\in \{1\ll n\}$.
For each such $j$, $(i,j)\notin Y'$ since $(i',j)\in F_0$; $(i,j)\notin X'$ since $(i,j)\notin F_0$; and similarly
$(i',j)\notin X', Y'$. By (3), $b_{ij}=a_{ij}$ and $b_{i'j} = a_{i'j}$, contradicting that $A$ is $t$-restricted.
From the symmetry between axes, this proves (4).
\\
\\
(5) {\em $B$ has breadth at most $4t$.}
\\
\\
Let $Z'=([m]\times[n])\setminus (X'\cup Y')$. Since $X\subseteq X'$ and $Y\subseteq Y'$, it follows that $Z'\subseteq Z$.
Suppose that for some $i,j$ all of the pairs $(i+h,j+h)\;(0\le h\le 4t)$ belong to $Z'$. Thus $(i,j)\in Z$, and
since $(i,j)\notin X'$, it follows from (3) that there exists $(i_1,j_1)\in P_1\cup Q_1$ such that
$i_1\le i$ and $j_1\le j$. Similarly, since $(i+4t,j+4t)\notin Y'$, there exists $(i_0,j_0)\in P_0\cup Q_0$
such that $i_0\ge i+4t$ and $j_0\ge j+4t$. But then $i_0\ge i_1+4t$ and $j_0\ge j_1+4t$, contrary to (2). This proves (5).
\\
\\
(6) {\em The local difference between $A,B$ is at most $4t$.}
\\
\\
Let $1\le i\le m$, and let $J$ be the set of all $j$ with $1\le j\le n$ such that $b_{ij}\ne a_{ij}$. We need to bound $|J|$.
Let $J_1$
be the set of all $j\in J$ such that $b_{ij} = 1$ and $a_{ij}=0$, 
and let $J_0$ be the set of all $j\in J$ such that $b_{ij} = 0$ and $a_{ij} = 1$.
If $j\in J_1$, then there is no $(i_0,j_0)\in P_0\cup Q_0$ such that $i_0\le i$ and $j_0\ge j$, and in particular
$(i,j)\notin Q_0$; and so $(i,j')\in E_0$ for at most $2t-1$ values of $j'< j$. Since $J_1\subseteq E_0$,
it follows (by choosing $j\in J_1$ maximal, if possible) that $|J_1|\le 2t$. 

If $j\in J_0$, then there is no $(i_1,j_1)\in P_1\cup Q_1$ such that
$i_1\le i$ and $j_1\le j$, and in particular $(i,j)\notin Q_1$, and so $(i,j')\in E_1$ for at most $2t-1$ values of $j'> j$;
and so $|J_0|\le 2t$. Summing, it follows that $|J|\le 4t$. From the symmetry between axes, this proves (6).

\bigskip
From (4), (5), (6), this completes the proof of \ref{matrix2}.~\bbox

\bigskip
Combining the previous two theorems, we deduce \ref{matrix}, which we restate::
\begin{thm}\label{matrixagain}
For all integers $t\ge 0$ and every $t$-restricted {\rm 0/1} matrix $A$,
there is an inclusive matrix $B$, such that the local difference between $A,B$
is at most $644t^4$.
\end{thm}
\Proof
By \ref{matrix2} there is a $t$-restricted {\rm 0/1} matrix $A'$ with breadth at most $4t$, such that the local difference between $A,A'$
is at most $4t$. By \ref{matrix1} applied to $A'$, with $w=4t$, there is an
inclusive matrix $B$ such that the local difference between $A',B$ is at most $640t^4$. But then the local difference between $A,B$ is at most
$640t^4+ 4t\le 644 t^4$. This proves \ref{matrixagain}.~\bbox

\section{Near-dominating graphs}

Now we use the results of the previous section to study $t$-dominating graphs.
It is tempting to try to apply \ref{matrix} directly to the adjacency matrix $A$ of a $t$-dominating graph, 
but that does not work. \ref{matrix} would give us an inclusive matrix $B$ that has bounded local difference from $A$, 
and so we could reorder its rows and columns to make it monotone; but we need the corresponding permutations of 
the rows and columns 
to be the same, and we need the monotone matrix to be symmetric, and the direct application of \ref{matrix} gives 
neither of these things. We will apply \ref{matrix} in another way.
First we show:

\begin{thm}\label{split}
Let $G$ be a $t$-dominating split graph. Then there is a split half-graph $H$ with $V(H)=V(G)$ such that the local difference between
$G,H$ is at most $644t^4$.
\end{thm}
\Proof
Let $G$ be a $t$-dominating split graph, and let $V(G)$ be the disjoint union of a clique $M$ and a stable set $N$.
Let $M=\{u_1\ll u_m\}$, where for $1\le i<j\le m$ the degree of $u_i$ is at most that of $u_j$, and let $N=\{v_1\ll v_n\}$, where
for $1\le i<j\le n$ the degree of $v_i$ is at most that of $v_j$. Let $A=(a_{ij})$ be the 0/1 matrix where $a_{ij} = 1$ if $u_i, v_j$
are adjacent and $a_{ij}=0$ otherwise. We claim that $A$ is $t$-restricted. For let $1\le i<i'\le m$, and suppose that there are
$t$ values of $j\in \{1\ll n\}$ such that $a_{ij} = 1$ and $a_{i'j} = 0$. Thus for each such $j$, $v_j$ is adjacent to $u_i$ and not
to $u_{i'}$. Since the degree of $u_{i'}$ is at least that of $u_i$, it follows that there are at least $t$ vertices adjacent to $u_{i'}$
and not to $u_i$, contradicting that $G$ is $t$-dominating. A similar argument shows that there do not exist $j,j'\in \{1\ll n\}$
with $j<j'$
such that for $t$ values of $i$, $a_{ij} = 1$ and $a_{ij'}=0$. Hence $A$ is $t$-restricted. By \ref{matrixagain}, 
there is an inclusive matrix $B$ such that the local difference between $A,B$
is at most $644t^4$. Let $H$ be the graph with vertex set $M\cup N$ in which $u_i, v_j$ are adjacent if $b_{ij}=1$, and $M$ is a clique
and $N$ is a stable set. It follows that $H$ is a split half-graph, and its local difference from $G$ is at most  $644t^4$.
This proves \ref{split}.~\bbox

\begin{thm}\label{reducetosplit}
Let $G$ be a $t$-dominating graph. Then there is a $t$-dominating split graph $H$ with $V(H)=V(G)$ 
and with local difference at most $2t$ from $G$.
\end{thm}
\Proof
Let $V(G)=\{v_1\ll v_n\}$, where for $1\le i<j\le n$, the degree of $v_j$ is at least that of $v_i$. Choose $i$ minimum such that
some vertex in $\{v_1\ll v_i\}$ has at least $2t+1$ neighbours in $\{v_1\ll v_i\}$. (If this is not possible then every vertex
has degree at most $2t$, and so $G$ has local difference at most $2t$ from the graph $H$ obtained by deleting all the edges, 
which is a $t$-dominating
split graph.) Let $N=\{v_1\ll v_{i-1}\}$ and 
$M = \{v_i\ll v_n\}$. Thus every vertex in $N$ has at most $2t$ neighbours in $N$. Choose $x\in N\cup \{v_i\}$ and 
$X\subseteq N\setminus \{x\}$ with 
$|X|\ge 2t$ (and $X\ne \emptyset$ if $t=0$), such that $X$ is the set of neighbours of $x$ 
in $N$. For all $j\ge i$, there are at 
most $t$ vertices adjacent to $x$ and not to $v_j$, since the degree of $v_j$ is at least that of $x$ and 
$G$ is $t$-dominating. In particular, $v_j$ is adjacent to at least half the vertices in $X$. Let $j\ge i$, and let $Y$ be the set of 
vertices in $M$ that are different from and nonadjacent to $v_j$. Since every vertex in $Y$ is adjacent to at least
half the vertices in $X$, and $X\ne \emptyset$, some vertex $v_h\in X$ is adjacent to at least half the vertices in $Y$. But there are
at most $t$ vertices adjacent to $v_h$ and not to $v_j$, since the degree of $v_h$ is at most that of $v_j$ and $G$ is 
$t$-dominating; and consequently $|Y|/2\le t$. Thus every vertex $v\in M$ is nonadjacent to at most $2t$
vertices in $M\setminus \{v\}$. Let $H$ be the split graph obtained from $G$ by deleting all edges with both ends in $N$
and making adjacent all pairs of vertices in $M$; then the local difference between $G,H$ is at most $2t$. 

We claim that $H$ is $t$-dominating. For let $u,v\in V(H)=V(G)$. If $u,v\in M$ then since one of $u,v$ $t$-dominates the other 
in $G$, the same is true in $H$ (since in $H$, $u,v$ have the same neighbours in $M\setminus \{u,v\}$). Similarly if
$u,v\in N$ then one of $u,v$ $t$-dominates the other. If $u\in M$ and $v\in N$ then $v$ $0$-dominates and hence $t$-dominates $u$.
This proves that $H$ is $t$-dominating, and so proves 
\ref{reducetosplit}.~\bbox

Combining \ref{split} and \ref{reducetosplit}, we obtain our main result, which we restate:
\begin{thm}\label{mainthm2}
Let $G$ be a $t$-dominating graph. Then there is a split half-graph $H$ with $V(H)=V(G)$ and with local difference at most 
$646t^4$ from $G$.
\end{thm}
\Proof
Let $G$ be $t$-dominating. By \ref{reducetosplit} there is a $t$-dominating split graph $G'$ with 
$V(G')=V(G)$ and with local difference at most $2t$ from $G$. 
By \ref{split}, there is a split half-graph $H$ with $V(H)=V(G)$ such that the local difference between
$G',H$ is at most $644t^4$. Thus the local difference between $G,H$ is at most $644t^4 + 2t\le 646t^4$.
This proves \ref{mainthm2}.~\bbox

\section{A counterexample}

Here is a concept similar to $t$-domination for graphs and $t$-restriction for matrices. Let $G$ be a bipartite 
graph with bipartition $(A,B)$; we say $G$ is
{\em $t$-nested} on $(A,B)$ if 
for all distinct $u,v\in B$, one of them has at most $t$ neighbours nonadjacent to the other.
Thus, $0$-nested graphs are half-graphs.
(Being $t$-nested is ``half'' of being $t$-restricted; $t$-restriction also requires the analogous statement 
with $A,B$ exchanged. We could also view $t$-nestedness as a version of $t$-domination for digraphs, if we regard the matrix as
the adjacency matrix of a digraph.) 

Let $G,H$ both be bipartite graphs with bipartition $(A,B)$; we say the {\em bipartite local difference}  between $G,H$ is 
the maximum over $v\in B$ of $|N_G(v)\setminus N_H(v)|+|N_H(v)\setminus N_G(v)|$.
One might hope that for all $t$ there exists $f(t)$ such that if $G$ is $t$-nested on $(A,B)$ then
there is a half-graph $H$ with the same bipartition $(A,B)$, with bipartite local difference at most $f(t)$ from $G$.
(Indeed, this looked like an easier question than our main result, and we tried it first as a warm-up.) But this is false,
even if $t=1$.
Here is a counterexample. 

Let $\mathcal{S}$ be the set of all  0/1 sequences with at most $k$ terms, and let $T$ be the binary tree defined by 
$\mathcal{S}$; thus $V(T) = \mathcal{S}$, and 
$s$, $s'$ are adjacent if their 
lengths differ by one and one of them is an initial subsequence of the other. Let $I$ be the set of all $s\in \mathcal{S}$
of length less than $k$, and $B$ those of length exactly $k$. (Thus $B$ is the set of leaves of $T$, and $I$ the set
of internal vertices.) 
Each $s\in B$ is a sequence of $k$ terms $(s_1\ll s_k)$ say, and we define $n(s)= \sum_{1 \le i\le k} s_i2^{k-i}$.
Thus all the numbers $n(s)\;(s\in B)$ are different and range from $0$ to $2^k-1$.
For $s,s'\in B$, we write $s<s'$ if $n(s)<n(s')$.

For each $s\in B$ take a set
$W_s$ of $2^{k+1}+2k$ new vertices, so that all the sets $V(T)$ and $W_s\;(s\in B)$ are pairwise disjoint.
Let $W= \bigcup_{s\in B} W_s$, and let $A=I\cup W$. We construct a bipartite graph $G$ with bipartition $(A,B)$ as follows.
Let $s\in B$. For $s'\in B$, $s$ has no neighbour in $W_{s'}$ if $s'>s$, and $s$ is adjacent to 
every vertex in $W_{s'}$ if $s'\le s$.
For $s'\in I$, we decide the adjacency of $s, s'$ in $G$ by the following rule. Choose $i$ minimum such that either
$i$ is greater than the number of terms of $s'$, or the $i$th terms of $s$ and of $s'$ are different, and let $x$ be the
$i$th term of $s$; we make $s,s'$ adjacent if and only if $x=0$.
(In other words, if we consider the paths $P_s$ and $P_{s'}$
from the root to $s$ and $s'$ in $T$, we have that $s,s'$ are adjacent in $G$ if the first vertex of $P_{s'}$ 
that does not lie on $P_s$ gets label 0.)

This graph is $1$-nested; and indeed, for $s,s'\in B$ with $s<s'$, there is at most one $v\in I$ adjacent to $s$ and not to $s'$,
and no such $v\in W$. Suppose that there is a $0$-nested graph $H$ with bipartition $(A,B)$ that has bipartite
local difference less than $k$ from $G$. For each $s\in S$, its degree in $G$ is at least $(2^{k+1}+2k)n(s)$ and 
at most $(2^{k+1}+2k)n(s) + 2^{k+1}-1$; and so its degree in $H$ is at least $(2^{k+1}+2k)n(s)-k$ and 
at most 
$$(2^{k+1}+2k)n(s) + 2^{k+1}-1+k< (2^{k+1}+2k)(n(s)+1)-k.$$
Consequently, for 
 distinct $s,s'\in B$, if $s<s'$ then the degree of $s$ in $H$ is less than that of $s'$. Since $H$ is $0$-nested,
it follows that every vertex of $A$ adjacent to $s$ in $H$ is also adjacent to $s'$.

Let $t_0$ be the null sequence, and let $L_0,R_0$ be the sets of members of $B$ with first term $0$ and $1$ respectively. 
It follows that $t_0$ is adjacent in $G$ to every vertex in $L_0$, and nonadjacent to every vertex in $R_0$; and yet $s<s'$ for every
$s\in L_0$ and $s'\in R_0$. Consequently, in $H$ either $t_0$ is nonadjacent to every member of $L_0$, or adjacent to every member of 
$R_0$. If the first, let $t_1$ be the one-term sequence $(0)$, and otherwise $t_1 = (1)$. Let $L_1$ be the set of members of $B$
such that $t_1$ is an initial segment, and the second term is $0$, and let $R_1$ be those such that $t_1$ is an initial segment
with second term $1$. Again, $t_1$ is adjacent in $G$ to every member of $L_1$, and nonadjacent to every member of $R_1$,
and in $H$ either $t_1$ is nonadjacent to all vertices in $L_1$ or adjacent to all in $R_1$; let $t_2$ be the corresponding 
two-term sequence. By continuing this process we obtain a sequence $t_0,t_1,t_2,\ldots,t_{k-1}\in I$, and a vertex $t_k\in B$,
where for $0\le i\le k$, $t_i$ has $i$ terms and $t_{i-1}$
is an initial segment of $t_i$, and for $0\le i\le k-1$, $t_k$ is adjacent to $t_i$ in exactly one of $G,H$. Since $t_k\in B$,
this contradicts that the bipartite local difference between $G,H$ is at most $k-1$.

This could be viewed another way: for each vertex in $B$, take its set of neighbours in $A$. Then we obtain a collection of
subsets $\mathcal{C}$ of $A$ such that for every two of them, say $X,Y$, one of $|X\setminus Y|,\allowbreak |Y\setminus X|\le~1$. 
But if we want to change this last $1$ to a $0$ by adding
and subtract elements of $A$ from the sets of $\mathcal{C}$,
then some set has to have
an arbitrarily large number of elements added or subtracted. 

\section{The Erd\H{o}s-Hajnal conjecture}

Let us say an {\em ideal} of graphs is a class $\mathcal{C}$ of graphs, such that if $G\in \mathcal{C}$ and $H$ is isomorphic
to an induced subgraph of $G$ then $H\in \mathcal{C}$; and an ideal is {\em proper} if some graph is not in it.
The Erd\H{o}s-Hajnal conjecture~\cite{EH} asserts:

\begin{thm}\label{EH}
{\bf Conjecture: } For every proper ideal $\mathcal{C}$, 
there exist $c,\epsilon>0$ such that for every graph $G\in \mathcal{C}$, $G$ has a clique or stable set of cardinality
at least $c|V(G)|^{\epsilon}$.
\end{thm}

We are interested in the way the (optimal) coefficient $\epsilon$ depends on $\mathcal{C}$. In particular, when does taking
$\epsilon = 1$ work? If $\mathcal{C}$ is the set of all graphs not containing one particular graph $H$ as an induced subgraph,
then there are almost no choices of $H$ for which $\epsilon = 1$ works --- only those graphs $H$ with at most two vertices, 
as is easily seen. 
For instance, let $S_t$
be the star with centre of degree $t$, and $\mathcal{C}$ the class of all graphs that do not contain $S_t$
as an induced subgraph; then $\epsilon = O(1/t)$, since the Ramsey number $R(k,t)$ is at least
$(k/\log k)^{(t+1)/2}$ for fixed $t$ and large $k$~\cite{keevash} (and $\epsilon$ is known to exist because of a result of
Alon, Pach and Solymosi~\cite{APS}).

But there are ideals defined by excluding more than one graph.  If we take $\mathcal{C}$ to be the class of graphs that contain
neither $S_t$ nor its complement as an induced subgraph, then~\cite{substars} for every graph $G\in \mathcal{C}$, 
either $G$ or its complement
has maximum degree bounded by a function of $t$, and so there exist $c,\epsilon$ with $\epsilon = 1$.

Another ideal of interest is the class of all $t$-dominating graphs, for fixed $t$. Every split graph has a clique or stable set containing at
least half its vertices, and so by \ref{reducetosplit}, in every $t$-dominating graph $G$,
there is a subset $X\subseteq V(G)$ with $|X|\ge |V(G)|/2$ such that either the subgraph induced on $X$ has maximum degree at most $2t$
or its complement graph does. Consequently $G$ has a clique or stable set with cardinality at least $|V(G)|/(4t+2)$;
and so we may take $\epsilon = 1$ in \ref{EH} for this class.

Take a ``substar'' (a graph obtained from a star by deleting some edges) and the complement of a substar,
and let $\mathcal{C}$ be the class containing neither of these graphs; in this case $\epsilon=1$ does not work in general (consider
a disjoint union of $n^{1/2}$ cliques each with $n^{1/2}$ vertices), but $\epsilon = 1/2$ works.
This is a special case of the following.
(If $G$ is a graph, $\omega(G)$ and $\alpha(G)$ denote the cardinalities of 
its largest clique and largest stable set respectively, and we denote $\max(\omega(G), \alpha(G))$ by $\rho(G)$.) 
We recall that threshold graphs were defined in section 2, and are the same as $0$-dominating graphs, because of~\ref{0dom2}.

\begin{thm}\label{thresholds2}
Let $H_1,H_2$ be threshold graphs, with $|V(H_1)|+|V(H_2)|=m$ and $\omega(H_1)+\alpha(H_2)=k$.
For every graph $G$, if $G$ has no induced subgraph
isomorphic to $H_1$ or to $H_2$, then $|V(G)|\le  (2^m-1) \rho(G)^{k-2}$.
\end{thm}
\Proof
We proceed by induction on $m$. If one of $H_1,H_2$ has at most one vertex the claim is trivial,
so we assume they both have at least two vertices.
A vertex is {\em isolated} if it has degree zero, and
{\em universal} if it is adjacent to all other vertices. For $i = 1,2$, let $H_i'$ be obtained from $H_i$ by deleting
an isolated vertex if there is one, and if not let $H_i' = H_i$; and let $H_i''$ be obtained from $H_i$ by
deleting a universal vertex if there is one, and if not let $H_i''=H_i$. 
Let $G$ contain neither of $H_1,H_2$.
For $X\subseteq V(G)$, $G[X]$ denotes the subgraph induced on $X$.
Let $v\in V(G)$, and let $N$ be the set of neighbours of $v$, and $M=V(G)\setminus (N\cup \{v\})$.

Suppose first that one of $H_1,H_2$ has an isolated vertex, and one has a universal vertex.
Since $G[M]$ contains neither of $H_1',H_2'$, it follows that $|M|\le (2^{m-1}-1) \rho(G)^{k-2}$;
and similarly $|N|\le (2^{m-1}-1) \rho(G)^{k-2}$. Consequently 
$$|V(G)|= |N|+|M|+1\le 2 (2^{m-1}-1) \rho(G)^{k-2} + 1\le (2^m-1) \rho(G)^{k-2}$$
as required. 

Thus either $H_1,H_2$ both have no isolated vertex, or they both have no universal vertex. We can reduce the second case to the first
by replacing
$H_1,H_2$ and $G$ by their complements, and exchanging $H_1,H_2$ (to preserve the property that $\omega(H_1)+\alpha(H_2)=k$); so
we may assume the first case holds.
It follows that $H_1,H_2$  both have
universal vertices, since they are threshold graphs;
and in particular, $\omega(H_1'')+\alpha(H_2'')= k-1$.

Choose $v\in V(G)$ with maximum degree, and let $N$ be its set of neighbours.
Since $G[N]$ contains neither of $H_1'', H_2''$, it follows from the inductive hypothesis that
$$|N|\le (2^{m-2}-1) \rho(G)^{k-3}\le 2^{m-2} \rho(G)^{k-3}-1.$$
Since $v$ has maximum degree, it follows that $G$ is $2^{m-2} \rho(G)^{k-3}$-colourable, and so some
stable set has cardinality at least 
$|V(G)|/(2^{m-2} \rho(G)^{k-3})$.
In particular $\rho(G)$ is at least this quantity, and so
$2^{m-2}\rho(G)^{k-2}\ge |V(G)|$, and again the result follows. This proves \ref{thresholds2}.~\bbox

\bigskip

One can also ask, what is the size of the largest clique or stable set
in  almost all graphs in an ideal? This may
be much larger than we can guarantee for {\em every} graph in the ideal.
Thus while for every $\epsilon>0$ there exists $t$ such that there are $K_t$-free 
$n$-vertex graphs with  no stable set of size at least $n^{\epsilon}$, it is known~\cite{EKR}
that for every
$t$ there exists $C_t>0$ such that almost every $K_t$-free $n$-vertex graph has a
stable set of size at least $C_t n$.

It remains an open problem to determine whether for some $\epsilon>0$ and
every graph $H$, there exists $C_H>0$ such that almost every graph $G$ that does not contain $H$ as an induced subgraph
has a stable set or clique of size at least $C_H |V(G)|^{\epsilon}$.
This would be true (see Lemma 3 of \cite{loebl} and the discussion
around it) if for some $\epsilon>0$ and every integer $t\ge 0$,
there
exists $C_t>0$ such that, for every graph $G$, either
\begin{itemize}
\item there are two vertices $u,v$ and a stable set $S$ with $|S\cap (N(u)-N(v))|$,
$|S\cap (N(v)-N(u))|$, $|S\cap N(u) \cap N(v)|$ all of cardinality $t$, or
\item there are two vertices $u,v$ and a clique $T$ with $|T\cap (N(u)-N(v))|$,
$|T\cap (N(v)-N(u))|$, $|T\setminus (N(u) \cup N(v))|$ all of cardinality $t$, or
\item $G$ contains a clique or stable set of size at least $C_t |V(G)|^\epsilon$.
\end{itemize}
Our results are a first step in proving this statement.


\begin{thebibliography}{99}
\bibitem{APS} N. Alon, J. Pach and J. Solymosi, ``Ramsey-type theorems with forbidden subgraphs'',
{\em Combinatorica} 21 (2001), 155--170.
\bibitem{keevash} T. Bohman and P. Keevash, ``The early evolution of the $H$-free process'', {\em Invent.
Math.} 181 (2010), 291--336.
\bibitem{substars} M. Chudnovsky, S. Norin, B. Reed and P. Seymour,
``Excluding a substar and an antisubstar'',
{\em SIAM J. Discrete Math.}, 29 (2015), 297--308.
\bibitem{EH} P. Erd\H{o}s and A. Hajnal, ``Ramsey-type theorems'', {\em Discrete Appl. Math.} 25 (1989), 37--52.
\bibitem{EKR} P. Erd\H{o}s, D. Kleitman and B. Rothschild, ``Asymptotic enumeration of $K_n$-free graphs'', 
{\em International Colloquium on Combinatorial Theory}, Atti dei Convegni Lincei 
17, Volume 2, Rome, 1976, 19--27.
\bibitem{FH} S. F\"oldes and P.L. Hammer, ``Split graphs'',
{\em Proceedings of the 8th South-Eastern Conference on Combinatorics, Graph Theory and Computing} (1977), 311--315.
\bibitem{loebl} M. Loebl, B. Reed, A. Scott, A. Thomason and S.
Thomass\'e,
``Almost all $H$-free graphs have the Erd\H{o}s-Hajnal property'', 
{\em An Irregular Mind  (Szemer\'edi is 70)}, in Bolyai Society Mathematical
Studies, Springer, Berlin, 21 (2010), 405-414.
\end{thebibliography}
\end{document}